\documentclass[11pt]{amsart}
 \usepackage{amsthm}
 \usepackage{amsmath}
 \usepackage{amsfonts}
 \usepackage{amssymb}
 \usepackage{amscd}
 \usepackage[all]{xy}
\usepackage[dvipsnames]{xcolor}
\usepackage{graphicx}
\usepackage{pdfsync}
\usepackage{url}

\usepackage{algorithm}
\usepackage{algorithmic}

\newtheorem{thm}{Theorem}[section]

\theoremstyle{definition}
\newtheorem{ex}[thm]{Example}

\newtheorem{prop}[thm]{Proposition}
\theoremstyle{definition}
\newtheorem{defn}[thm]{Definition}

\newtheorem{example}[thm]{Example}
\numberwithin{equation}{section}

\def \OF { \Psi }
\def \OFu { \Psi_u }

\def \Xreg {X_{reg}}

\def \JacX {\text{Jac}X}

\DeclareMathOperator{\minors}{minors}

\def \cL { {\mathcal L} }

\def \A { {\mathbb A} }

\def \R { {\mathbb R} }
\def \C { {\mathbb C} }

\def \P { {\mathbb P} }

\newcommand{\<}{\langle}
\renewcommand{\>}{\rangle}

\DeclareMathOperator{\rank}{rank}

\newcommand{\lhra}{\ensuremath{\lhook\joinrel\relbar\joinrel\relbar\joinrel\rightarrow}}

\newcommand{\demph}[1]{\emph{\color{Blue}#1}}

\newcommand{\abr}[1]{ {\color{Orange}{ #1  --Abraham}}} 

 \title{Critical points via monodromy and local methods}

\author[Mart\'in del Campo]{Abraham Mart\'in del Campo}
\address{
Abraham Mart\'in del Campo\\
IST Austria\\
Am Campus 1\\
3400 Klosterneuburg, Austria}
\email{abraham.mc@ist.ac.at}
\urladdr{http://pub.ist.ac.at/~adelcampo}

\author[Rodriguez]{Jose Israel Rodriguez}
\address{
Jose Israel Rodriguez\\
Department of Department of Applied and Computational Mathematics and Statistics\\
University of Notre Dame\\
255 Hurley \\
Notre Dame, IN 46556 \\
USA
}
\email{jo.ro@ND.edu}
\urladdr{http://www.nd.edu/~jrodri18/}

\thanks{This research paper is partly supported by 2014 NIMS Thematic Program on Applied Algebraic Geometry.}
\thanks{The second author is partially supported by the National
Science Foundation under Award No. DMS-1402545.}

\begin{document}

\begin{abstract}
In many areas of applied mathematics and statistics, it is a fundamental
problem to find the best representative of a model by
optimizing an objective function. 
This can be done by determining critical points of the objective function restricted to the model. 

We compile ideas arising from numerical algebraic geometry 
to compute the critical points of an objective function.
Our method consists of using numerical homotopy continuation 
and a monodromy action on the total critical space to compute 
all of the complex critical points of an objective function. 
To illustrate  the relevance of our method, we apply it
to the Euclidean distance function to compute ED-degrees and  the likelihood function to compute   maximum likelihood degrees.
%
\end{abstract}

\maketitle

\section{Introduction}

In science and engineering, it is common to work with models that can be described as the solutions of a parametrized  system of polynomial equations. 
For such algebraic models $X\subset \R^n$, we are interested in the following polynomial optimization problem: given $u\in \R^n$, find those points $x^*\in X$ that optimize an objective function $\OF(x,u):\R^n \times \R^n \to \R$.  Common objective functions seen in applications include norms, distances, and other statistical functions. In this paper, we focus on the case when $\OF$ is a quadratic norm (such as the  Euclidean distance) and when $\OF$ is the likelihood function.

By Fermat's theorem, we find $x^*$ by computing the set of critical points using the method of Lagrange multipliers.  Thus, we are interested in finding the points $x\in X$ for which the gradient $\nabla\OF(x,u)$ is orthogonal to the tangent space $T_xX$ of $X $ at the point $x$. 
We recast these constraints as the solutions to a system of polynomial equations.
In the algebraic closure, the number $d$ of critical points is an invariant that gives information about the algebraic complexity of the optimization problem. Finding the number $d$ for different objective functions and models is 
a challenging problem. When $\OF$ is the Euclidean distance, the number $d$ is called the Euclidean distance degree (ED-degree), and when $\OF$ is the likelihood function, the number $d$ is the maximum likelihood degree (ML-degree).

Finding these degrees $d$ is an active area of research. 
There are already some fundamental results about these algebraic degrees and the geometry behind them.
A general degree theory for the ED-degree was introduced in ~\cite{DHOST13} and for the ML-degree in~\cite{CHKS06,HKS05, HS13}.
For some special cases, it is possible to find formulas for $d$ using techniques from algebraic geometry~\cite{AAGL14, Huh13,HS13, Lee14, NR09, OttSpaStu14,Pie15,Rod14}, some are summarized in~\cite{Ran12}. However, this is not always the case, and other methods are necessary. 

There are  algorithms proposed for identifying critical points based on local methods. Such methods involving Gr\"obner bases include ~\cite{FauSafSpa12, OttSpaStu14, Spa14}. From the numerical algebraic geometry perspective, those in~\cite{WR13, BBHSW14} are can be used and involve  SVD decomposition and regeneration (respectively) to find critical points for the computation of witness sets. 
For computing the real critical points, one method is to compute all complex solutions and determine the real solutions among them. 
This can be avoided by using local methods to determine local optima and using  heuristics to decide if the local optima is in fact global.
In these cases, one can study the expected number of real critical points when $u$ is drawn from a given probability distribution.
This has been 
 studied for the ED-degree under the name of \emph{average ED-degree} (e.g.~\cite{DLR15}).
Recently, a probabilistic method was proposed in~\cite{RT15}, where they classify the real critical points for the likelihood function.

In this paper, we propose a different numerical method to compute all the complex critical points.
Our method consists of two parts and complements those used in~\cite{GR13, HRS12}.  
The first part consists of exploiting a monodromy group action to randomly explore the variety of critical points. 
In this random exploration, we determine a subset of  critical points. 
The second part consists of a trace test to determine, with probability one, if the computed subset contains \emph{all} of the critical points. 

This paper is organized as follows.
In the next section, we give a geometric formulation of the critical equations along with a concrete formulation of the critical equations. 
In the following section, we review monodromy for a parameterized polynomial system.  
This includes a description of each part of our method, including a trace test. 
We conclude by using these methods to reproduce some known  ED-degrees and ML-degrees, and compare the times of those computations against our method.
%
Throughout the paper, we include  implementation subsections so the reader can use the methods in their own research. Supplementary materials can be found on the second authors website at \texttt{http://www3.nd.edu/}$\sim$\texttt{jrodri18/monodromy}.

%
%
We end this introduction with an illustrating example to  the critical points problem.
\begin{ex}\label{firstExample}
Let $X$ be the ellipse defined as the algebraic variety of the points $x=(x_1,x_2)\in \R^2$ that satisfy the polynomial equation
\begin{equation*}\label{eq:Felipse}
1744x_1^2-2016x_1x_2-2800x_1+1156x_2^2+2100x_2+1125=0.
\end{equation*}
For a (generic) choice of $u=(u_1,u_2)\in \R^2$, we are interested in the optimization problem
\begin{equation}\label{eq:optimizedistance}
\min \OF(x,u)=(x_1-u_1)^2+(x_2-u_2)^2 \quad \mbox{subject to }x\in X.
\end{equation}
Here, the objective function $\OF(x,u)$ is the square of the Euclidean distance between the (fixed) point $u\in \R^2$ and the ellipse $X$. 
The critical points of~\eqref{eq:optimizedistance} are those points $x^*\in X$ whose tangent is perpendicular to the line segment joining $x^*$ and $u$. 
This is illustrated in Figure~\ref{F:elipse}. 
\begin{figure}[htb]
 \includegraphics[height=0.31\textheight]{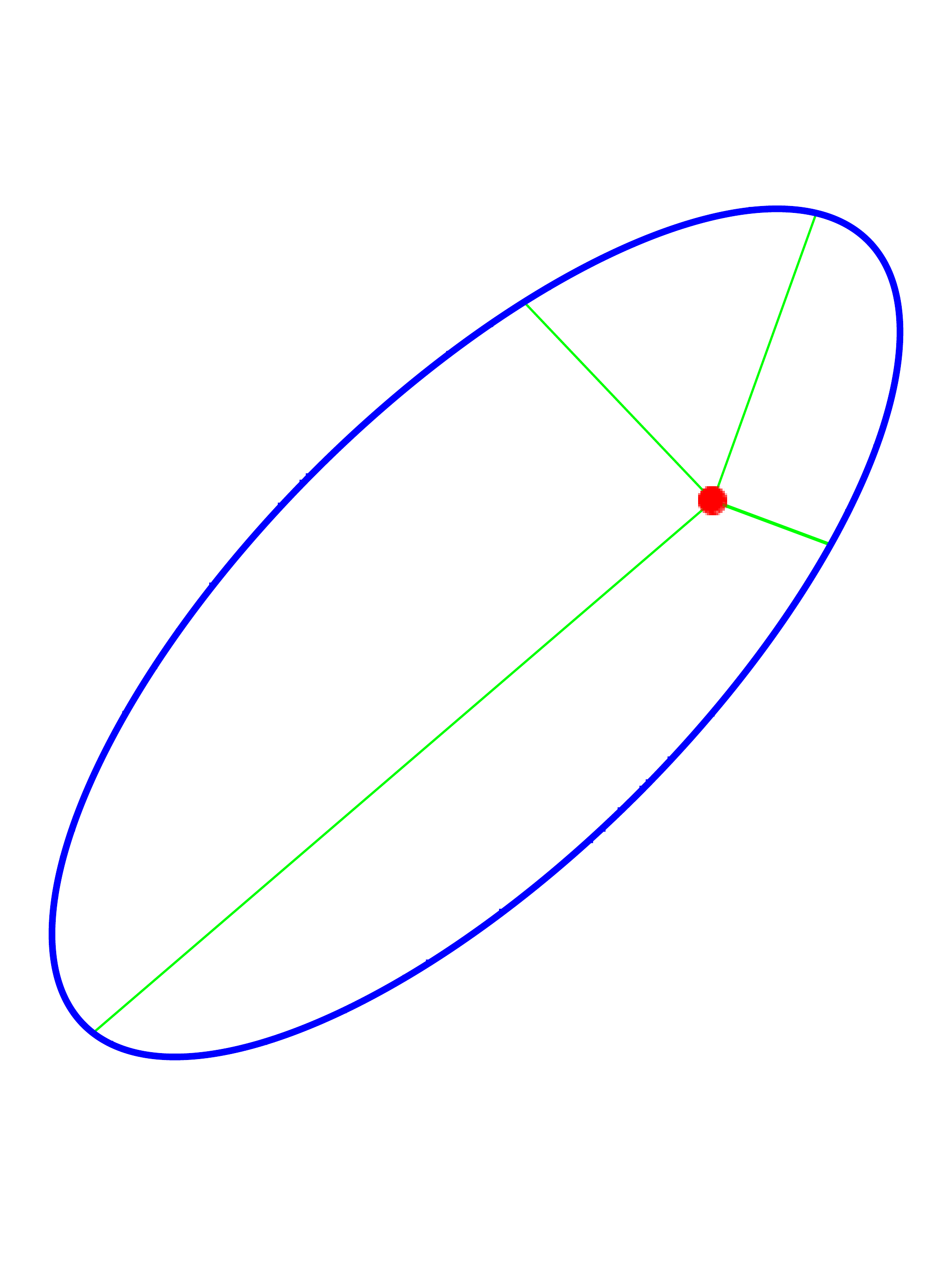}
 \vspace{-15pt}
 \caption{Critical points for the distance between  the point $u=(0.75, -0.29)$ and the ellipse  $1744x_1^2-2016x_1x_2-2800x_1+1156x_2^2+2100x_2+1125=0 $ }
 \label{F:elipse}
\end{figure}
By using Lagrange multipliers, finding the critical points corresponds to finding the solutions to the following polynomial system 
\begin{eqnarray}\label{eq:criticalsystem}
\nonumber 1744x_1^2-2016x_1x_2-2800x_1+1156x_2^2+2100x_2+1125 &=&0\\
x_1-u_1 + \lambda \left( 3488 x_1-2016 x_2 - 2800 \right) &=& 0\\
\nonumber x_2-u_2 + \lambda \left(-2016 x_1+ 2312 x_2 + 2100 \right) &=&0.
\end{eqnarray}
where $x_1, x_2,$ and $\lambda$ are indeterminates, and $u_1, u_2$ are parameters.
For instance, when $u=(0.75, -0.29)$, we use tools from numerical algebraic geometry to find the 
set of critical points, which consist of the points 
$$(0.8444, -0.3330), \, (0.8329, -0.0066), \, (0.5985, -0.0941), \, \mbox{and }\, (0.2529, -0.8140).$$ 
These are the coordinates on $x_1$ and $x_2$ of the solutions to \eqref{eq:criticalsystem}. In this case, there are only four solutions, showing that the ED-degree for $X$ is 4. 
The solution to the optimization problem~\eqref{eq:optimizedistance} corresponds to the critical point with the smallest Euclidean distance to $u$, which in this case is the first point. 
\end{ex}

\section{Critical points}\label{s:criticalpoints}
In this section, we will introduce critical points of $\OF$ restricted to an affine variety  $X$  defined by a system of polynomials in the indeterminates 
$x_1,\dots,x_n$.
To a given variety $X\in \C^n$ and a function $\OF$, we associate its \emph{total critical variety} $Y \in \C^n\times \C^n$ that consists of all the critical points of $\OF$ in $X$. This critical variety is defined by a system of polynomial equations in two sets of indeterminates: $x_1,\ldots, x_n$ and $u_1,\dots,u_n$. 
We refer to $x_1,\dots,x_n$ as the indeterminates of the model  and to $u_1,\dots,u_n$ as the parameters for $\OF_u$.

We let the system of polynomials $F=\{f_1,\dots,f_m\}$ in the indeterminates $x_1,\dots,x_n$ define an affine variety $X$ of $\C^n$. 
We call this variety the \demph{model}. 
In engineering and statistics, we are interested in the real points of the variety $X$. 
However, we work over the complex numbers and then restrict to the real numbers, 
 as is usual in applied algebraic geometry.
  In this way, the monodromy methods that we discuss in the next sections will let us conclude properties about $X$ that can be used to understand its real counterpart.

The objective functions $\Psi$ we  consider   are the Euclidean norm and the  likelihood function, but our method holds also for functions $\OF(x,u): \R^n\times \R^n \to \R$ with a derivative  (or logarithmic derivative) that is a rational function in $x_1,\dots,x_n$.
We regard $\OF$ as a family of functions $\OF_u$ parametrized by $u\in \C^n$.


Let $\JacX$ denote the \demph{Jacobian of $X$}, 
which is the $n{\times} m$-matrix whose $(i,j)$th entry is given by $\partial f_j/\partial x_i$. 
Let $\nabla f_j$ denote the \demph{gradient vector of} $f_j$; thus, we can write
$$
\JacX =
\left[\begin{array}{ccc}
\nabla f_{1}^\top & \cdots & \nabla f_{m}^\top
\end{array}\right].
$$
The rank of the Jacobian  of $X$ evaluated at a point $x$ in $X$ is at most the codimension of $X$.   
We say the point $x^*$ in $X$ is \demph{regular} if the Jacobian at $x^{*}$ has rank equal to the codimension of $X$; 
otherwise, the point $x^*$ is said to be \demph{singular}. 
We let $\Xreg$ denote the set of regular points of $X$ and
let  $X_{sing}$ denote the set of singular points of $X$. 
Thus, we have $X = \Xreg \sqcup X_{sing}$.

If $\nabla \OF$ denotes the (truncated) gradient vector whose entries are the partial derivatives $\partial \OF / \partial x_i$ for all $i=1,\ldots,n$, then, a critical point of $\OF$ on $X$ can be defined with respect to the Jacobian of $X$ and $\nabla\OF$.
Let $H$ be the set of points in $\C^n$ where any denominator of the rational coordinates of $\nabla\OF$ vanishes. 
\begin{defn}\label{d:critical_point}
A point $x \in X$ is a \demph{critical point of } $\OF$ if and only if 
\begin{equation}\label{eq:crit_conditions}
x \in \Xreg\backslash H \quad \mbox{and }\quad 
\rank \left[\begin{array}{cccc}
\nabla \OF^\top &
\nabla f_{1}^\top &
\cdots &
\nabla f_{m}^\top
\end{array}\right]
\leq  m.
\end{equation}
\end{defn}

When $\OF_u$ is the likelihood function 
$\OF(x,u) = x_1^{u_1}\cdots x_n^{u_n}$, the gradient  $\nabla\OF(x,u)$ equals $(\frac{u_1}{x_1},\frac{u_2}{x_2},\dots,\frac{u_n}{x_n})$; thus, the associated $H$ is the solution set defined by  $x_1x_2\cdots x_n=0$.
For the quadratic norm $\OF(x,u) = (x_1-u_1)^2+(x_2-u_2)^2+\cdots +(x_n-u_n)^2$, 
the gradient is $\nabla\OF(x,u)=(x_1-u_1,x_2-u_2,\dots,x_n-u_n)$; thus, $H$ is the empty set.

Our problem consists of finding those $x^*$ such that,
\begin{equation}\label{eq:optimization}
\min / \max \ \OF_u(x) \quad \text{subject to } x\in X\backslash H, \text{ for fixed }u\in \R^n.
\end{equation}
%
%
Note that because the conditions~\eqref{eq:crit_conditions} are additive with respect to irreducible components we can assume that $X$ is an irreducible variety.
We regard the objective function $\OF$ as a function depending on a fixed parameter $u\in \R^n$; thus, we interpret it as a family of objective functions $\OFu$ parametrized by $u$. Therefore, we are interested in the set of critical points for the parametrized family of objective functions $\OFu$. 
\begin{defn}\label{d:critical_variety}
Let $X$ be an irreducible variety and $\OF_u$ be a parametrized family of objective functions as above. We define the \demph{total critical variety $Y$} as the Zariski closure of the set of critical points of $\OF_u$ in $\Xreg\backslash H$, this is 
\begin{equation}\label{eq:crit_variety}
Y := \overline{\{ (x,u)\in \C^n \times \C^n : x\in \Xreg\backslash H \mbox{ and } x \mbox{ is a critical point for }\OF_u \}}.
\end{equation}
\end{defn}
The total critical variety has a natural projection to $\C^n$ associated to the coordinates $u_1,\dots,u_n$. 
When $\OF_u$  is the likelihood function, $Y$ is called the \textit{likelihood correspondence}, and the degree of this projection over a general point is called the \textit{maximum likelihood degree (ML degree) of}  $X$.  
When $\OF_u$ is the Euclidean distance function, $Y$ is called the $ED correspondence$, and the degree of this projection over a general point is called the \textit{Euclidean distance degree (ED degree) of} $X$.

The constraints \eqref{eq:crit_conditions} impose algebraic equations that exhibit $Y$ as an algebraic variety. Working in the polynomial ring $\C[x_1,\ldots, x_n,u_1,\ldots, u_n]$, we let $\mathcal I(X)$ denote the prime ideal defining $X$,
let $\mathcal I(H)$ denote the ideal of $H$, 
and let $\mathcal J(X,\OF)$ denote the $(m{+}1)$-minors of the matrix 
\begin{equation*}\label{eq:matrix_ext_jac}
\left[\begin{array}{cccc}
\nabla \OF^\top &
\nabla f_{1}^\top &
\cdots &
\nabla f_{m}^\top
\end{array}\right].
\end{equation*}
We regard $\nabla \OF$ as the vector whose entries are only the partial derivatives $\partial \OF/\partial x_i$ for $i=1,\ldots, n$.
Thus, the defining equations for $Y$ are the generators of the ideal
\begin{equation}\label{eq:overdetermined_equations}
( \mathcal I(X) + \mathcal J(X,\OF)) : (\mathcal I (H)\cdot\< \minors(m, \JacX)\>)^\infty.
\end{equation}
These equations form an overdetermined system of equations. However, since our techniques use numerical algebraic geometry tools, we would prefer to define the total critical variety by a square system of equations where the number of 
indeterminants  equals the number of equations. 

One way to do this is to  introduce auxiliary unknowns and consider a variety 
$\widehat Y$ that is an irreducible component of a variety $\widehat Z$ where $\widehat Y$ projects to $Y$.
We use Lagrange multipliers to define a squared systems as follows. For a fixed $u\in \C^n$, consider the polynomial system $G:\C^n\times \P^m \to \C^{m+n}$ given by 
\begin{equation}\label{eq:local_system}
G(x,\lambda) := \left[\begin{array}{c} 
F(x) \\
\lambda_0 \nabla \OF_u(x)^\top + \lambda_1 \nabla f_1(x)^\top +\cdots + \lambda_m \nabla f_m(x)^\top
\end{array}\right].
\end{equation}
Note that if $x^*\in X$ satisfy the conditions~\eqref{eq:crit_conditions}, then there exists $\lambda^*\in \P^m$ such that $G(x^*, \lambda^*)=0$ by the Fritz John condition~\cite{John48}. In the affine patch where $\lambda_0 =1$, the system $G$ becomes a square system and its solutions $(x^*,\lambda^*)$ project to  critical points $x^*\in \Xreg$. 
We will use the system~\eqref{eq:local_system} when we refer to the square system defining the total critical variety $Y$.

\begin{example}\label{tcExample}
Consider $X$ in $\C^4$ defined by the equations 
$f_1=x_1x_3-x_2^2,f_2=x_2x_4-x_3^2,f_3=x_1x_4-x_2x_3$ and the objective function $\OF_u=(x_1-u_1)^2+\cdots+(x_4-u_4)^2$. 
The codimension of $X$ is $2$. 
Let $F=\{\widehat f_1, \widehat f_2\}$ be the following linear combinations of the $f_j$:
$$
\widehat f_1  = \frac{1}{100}   (2f_1 +3f_2+5f_3) \qquad 
\widehat f_2  = \frac{1}{100}   (7f_1 +11f_2+13f_3)
$$
The variety defined by $F$ is reducible but still contains $X$ as an irreducible component. 
The additional component $X'$ is the 2-dimensional linear space defined by 
$9x_2-x_3-16x_4=81x_1-145x_3-16x_4=0$.
This randomization procedure of making a regular sequence from an overdetermined system of equations is a standard tool of numerical algebraic geometry. 
The equations from $(\ref{eq:local_system})$ define a reducible variety $\widehat Z$ 
in $\C^2\times \C^4\times \C^4$. 
This is a system of $6$ equations in $10$ indeterminants $x_1,x_2,x_3,x_4,\lambda_1,\lambda_2,u_1,u_2,u_3,u_4$ given by
\[
\begin{array}{cccc}
\widehat f_1=0,\quad \widehat f_2 =0\\
(x_1-u_1)+\lambda_1\frac{\partial \widehat f_1}{\partial x_1}+\lambda_2\frac{\partial \widehat f_2}{\partial x_1}=0\\
\vdots \\
(x_4-u_4)+\lambda_1\frac{\partial \widehat f_1}{\partial x_4}+\lambda_2\frac{\partial \widehat f_2}{\partial x_4}=0.
\end{array}
\]
This variety is reducible with two components $\widehat Y$ and $\widehat Y'$.
The component $\widehat Y$  and the total critical variety $Y$ of $X$ are birationally equivalent, and the 
same is true for 
 the second component $\widehat Y'$ and the total critical variety $Y'$ of the linear space $X'$. 

The fiber of the projection of $\widehat Y$ and thus $Y$ to $\C^4$ associated to the coordinates $u_1,u_2,u_3,u_4$ over the point $u^*=(\frac{2}{5},-\frac{2}{7},\frac{5}{6},\frac{3}{7})$ is described by the $6$ equations above and the $4$ equations~below:
$$u_1=2/5, \quad u_2=-2/7,\quad u_3=5/6,\quad u_4=3/7.$$
Solving the system of $10$ equations, we find $7$ solutions. 
Six of them have $x$ coordinates in $X$ and the last solution has $x$ coordinates in $X'$.
For this example, there are only 3 real critical points and their $x$ coordinates are listed below. 
The first two points lie in $X$ while the last point is in  $X'$.
\[
\begin{array}{ccccc} 
x_1&x_2&x_3&x_4& \Psi_u(x)\\
.128515& .252579& .496407& .  975616& .776241\\
.365062& .0690207& .0130495& .00246721& .981488 \\
.651048& -.288682& .384257&   -.186399&.642893.
 \end{array}
 \]
\end{example}
The first and last points in the list are the closest points from $u^*$ in $X$ and $X'$ respectively. 
The following theorem justifies that the ML-degree and ED-degree are well defined. The corresponding proofs can be found in~\cite{DR14} and~\cite{DHOST13} respectively.

\begin{thm}\label{thm:totvarietyirred}
The total critical variety $Y$ is an irreducible variety of dimension $n$ inside $\C^n\times \C^n$ and there exist open sets where the second projection $p:Y\to \C^n$ is generically finite and dominant.
\end{thm}

With our monodromy techniques, we were able to compute the $6$ points in the fiber of $u^*$ from $p: Y\to \C^4$ in Example \ref{tcExample}.  These techniques have the advantage of ignoring the critical point on \textit{junk  components} (in this case $X'$).  
The novelty of this paper is the compilation of ideas arising from numeral algebraic geometry applied to the computation of critical points.
In addition, we have developed an implementation of monodromy homotopies that can be used to compute fibers 
of  projections. 
This implementation uses the numerical algebraic geometry software \texttt{Bertini}~\cite{Bertini}  and commutative algebra software \texttt{Macaulay2}~\cite{M2}. In particular, it involves the software packages \texttt{Bertini.M2}~\cite{Bertini4M2} and \texttt{NAGtypes}~\cite{NAG4M2}.

\section{Monodromy and general methods}

In this section we define the action of a monodromy group on the fiber of the projection of the total critical variety $Y$. This action is one of the  main tools  we use to compute critical points. 
We start by considering a polynomial system $F=\{f_1,\dots,f_m\}$ that defines an irreducible affine variety $X$ of codimension $k$ and an objective function $\OF$.
The variety $X$ and the function $\OF$ 
determine the total critical variety $Y$ from Definition~\ref{d:critical_variety}. 
The projection $\pi: Y\to X\subset \C^n$ given by $(x,u)\to x$ realizes $Y$ as an affine vector bundle of rank $k$ over $\Xreg$.  
For the second projection $p:Y \to \C^n$ defined by $(x,u)\mapsto u$, there is an open set $U\subset \C^n$ where the following fiber diagram holds
 \begin{equation}\label{Eq:fiber_diagram}
  \raisebox{-20pt}{
  \begin{picture}(60,45)
   \put(-7,35){$p^{-1}(U)$} \put(28,35){$\lhra$} \put(55,35){$Y$}
   \put(-3,19){$p$}\put(9,32){\vector(0,-1){20}}
      \put(58,32){\vector(0,-1){20}}\put(62,19){$p$}
   \put(5, 0){$U$} \put(22, 0){$\lhra$} \put(54, 0){$\C^n$}
  \end{picture}
  }
 \end{equation}
with the restriction of $p$ onto $U$ being a generically finite morphism of degree $d$. 
A loop in $U$ based at $u$ has $d$ lifts to $p^{-1}(U)$, one for each point in the fiber $p^{-1}(u)$. 
Associating a point in the fiber $p^{-1}(u)$ to the endpoint of the
corresponding lift gives a permutation in $S_d$. This defines the usual permutation action
of the fundamental group of $U$ on the fiber $p^{-1}(u)$. The \demph{monodromy group} of the map
$p:Y \to \C^n$ is the image of the fundamental group of $U$ in $S_d$.

The equations of (\ref{eq:local_system}) define a reducible variety $\widehat Z$ which contains 
an irreducible component $\widehat Y$ which is birationally equivalent to $Y$.
Both $\widehat Z$ and $\widehat Y$ lie on $\C^m\times\C^n\times \C^n$, and the projection 
$\hat p:\widehat Y \to \C^n$ to the $u$ coordinates factors with the projection $p: Y \to \C^n$.
This factorization is compatible with the monodromy action on the fiber $\hat p^{-1}(u)$ over a regular point $u$.

The idea behind the method is to use numerical homotopy continuation to compute the fiber at a regular point $u^*\in U$. 
Suppose that we are endowed with a point  $(\lambda^*_0,x^*_0,u^*_0)$ in the fiber $\hat p^{-1}(u^*_0)$. 
We generate a random loop $\gamma:[0,1] \to U$ based at $u$, meaning $\gamma(0) = \gamma(1) = u^*_0$. 
 We numerically follow the points  in the fibers $\hat p^{-1}(\gamma(t))$ as $t$ deforms from 0 to 1. 
 This computes a lift of $\gamma$ to  a path from $x_0^*$ to another point $x^*_1$ in the fiber $\hat p^{-1}(u^*_0)$. Computing sufficiently many of these random loops enable us to recover the fiber.  

\subsection{Populating the fiber}\label{ss:pop_fiber}

To describe our method, we start by discussing its components and the numerical tools they use. 
Suppose for a general  $u$ we have a critical point $x_0$ for $\OF_u$. 
Then, we create a random loop $\gamma$ based on $u$ in the parameter space, and we use parameter homotopies to track the path from the point $x_0$ to another point $x_1$ in the fiber $p^{-1}(u)$. We now describe the details of this part.

Let $\phi(x; u):\C^N\times \C^k \to \C^N$ be a parametrized family of polynomial equations
$$
\phi(x,u) := \left[\begin{array}{c} 
\phi_1(x_1, \ldots, x_N; u_1,\ldots , u_k) \\
\vdots \\
\phi_N(x_1, \ldots, x_N; u_1,\ldots , u_k)
\end{array}\right] = 0.
$$
In other words, $\phi(x;u)$ is a system of $N$ polynomials in the variables $x\in \C^N$ and parameters $u \in \C^k$. 
Thus, for a fixed choice of $u\in \C^k$, we have a squared system. The number of nonsingular (isolated) solutions of $\phi(x,u)=0$ remains constant for general choices of parameter values $u\in \C^k$. 
Parameter homotopy consists of tracking
known solutions of the system $\phi(x,u)$ for specific parameter values $u\in \C^k$ to find solutions of the system for other parameter values $u' \in \C^k$. For a fixed $u\in \C^k$, a  path between $u$ and $u'$ is a continuous function $\gamma(t): [0,1] \to \C^{k}$  such that $\gamma(0) = u$ and $\gamma(1)=u'$, for which $\gamma(t)\in \C^k$ stays generic for $t\in [0,1)$. If $S$ is the set of solutions of $\phi(x;u)$ for a generic choice of parameters $u\in \C^k$, a \demph{parameter homotopy} consists of choosing a path $\gamma(t)$ between $u$ and a generic  $u'$, and follow the solutions of the homotopy
$$
h(x,t) := \phi(x,\gamma(t)) = 0
$$
starting at $S$ as $t$ goes from 0 to 1. In particular, if $\gamma$ is a loop based at $u$, then $h(x,t)$ defines a path between the solution set $S$, which we call a \demph{monodromy path}. 

We refer the reader to~\cite{BHSW13, SW05} for more about parameter homotopies and we focus on its application to our problem.
 The parametrized system $\phi$ that we consider is the one defined in the affine chart where $\lambda_0=1$ by $G(x,\lambda; u)=0$ as in~\eqref{eq:local_system}:
$$
G(x,\lambda; u) := \left[\begin{array}{c} 
F(x) \\
\nabla \OF(x;u)^\top + \lambda_1 \nabla F_1(x)^\top +\cdots + \lambda_m \nabla F_m(x)^\top
\end{array}\right],
$$
a system of $n{+}m$ polynomials in the variables $x_1,\ldots, x_n, \lambda_1,\ldots ,\lambda_m$ and parameters $u_1,\ldots, u_n$.
We assume that for a given choice of parameters $u\in \C^n$, we know a critical point $x^*\in \C^n$. 
In Section~\ref{ss:1st_critpt}, we discuss some possible ways to find the first critical point $x^*$ for a given choice of parameter values $u\in \C^n$; for now, we focus on the way we use parameter homotopies to find the rest of the critical points. 

Let $F:\C^n \to \C^m$ be a polynomial system defining our model and let $S$ denote its complete set of critical points. By letting $\lambda^* = (0,\ldots, 0)\in \C^m$, the system $G$ vanishes at $(x^*,\lambda^*)$, so
let $S_0 := \{ (x^*,\lambda^*)\}$ be the starting solution set. We generate a loop based on $u$ by taking two random parameter points $u', u''\in \C^n$ and we create a triangular loop $u \to u' \to u'' \to u$ using linear paths of the form $t\cdot u'' + (1-t)\cdot u'$.
We track the solution $(x^*,\lambda^*)$ to a new solution $(x', \lambda')$ of $G(x,\lambda; u)$.
In this way, if $S_r\subseteq S$ is a set of solutions obtained by repeating this process $r$ times, we 
obtain a new solution set by taking a new random loop to track $S_r$ to $S'_r$ and letting $S_{r+1} := S_r \cup S'_r$. Notice that $S'_r$ may coincide with $S_r$, but we always have the inclusions $S_r \subseteq S_{r+1}\subseteq S$. We continue constructing random monodromy paths until $S_{r+1} = S$. 
To verify this last step, we perform a trace test, which we explain next.

\subsection{Trace test}\label{ss:trace_test}

The trace test was introduced in~\cite{SVW2002} where they use it in a method for computing the solutions of a polynomial system using monodromy and to decompose positive dimensional reducible varieties.
Here, we state the trace test criterion and we illustrate it with an example. 

Suppose we have a reduced irreducible  $1$-dimensional affine variety $X$ of $\C^n$.
The intersection of $X$ with a generic hyperplane $\cL$ consists of $\deg(X)$-many points. The trace of X with respect to the hyperplane $\cL$ is the point defined by the coordinate-wise sum of the points in $\cL\cap X$.
If $\cL$ is defined by a linear form $l(x)$, let $\cL_t$ be a linear deformation of $\cL$ induced by $l(x)+t$.
The trace of $\cL_t$  depends on $t$. The main result of \cite{SVW2002} says the following. 

\begin{prop}
With the assumptions and notation above, the trace of $X$ with respect to the family of hyperplanes $\cL_t$ is affine linear in $t$.  
Moreover, the coordinate-wise sum of any non-empty proper subset of 
$\cL_t\cap X$ is not affine linear in $t$. 
\end{prop}

The idea behind the trace test is that if $X$ is not linear, a linear deformation $\cL_t$ will not deform the points $\cL_t\cap X$ linearly, but their trace will.
In monodromy methods, the trace test is particularly useful to verify that all points of $\cL\cap X$ are found. 
We illustrate this test with the following example.

\begin{example}
Consider $X$ defined by $x_1^2-x_2=0$.
Let $\cL$ be defined by $l(x)=2x_1+4x_2-1$. 
The trace of $X$ with respect to the linear deformation $\cL_t$ induced by $l(x)+t$ is 
$( -\frac{1}{2}, -\frac{1}{2}t+\frac{3}{8} )$.
However, the nonempty proper subsets of $\cL_t\cap X$ have  traces equaling the points 
$$
\left(\tfrac{-1+\sqrt{-4t+5}}{4},
\tfrac{-(2t-3)-\sqrt{- 4t + 5  }}{8}\right)
\text{ and }
\left(\tfrac{-1-\sqrt{-4t+5}}{4},
\tfrac{-(2t-3)+\sqrt{- 4t + 5  }}{8}\right).
$$
\end{example}

The results of the forthcoming work \cite{HR15} give a trace test for 
multi-projective varieties. Dehomogenizing the multi-projective variety allows us to give a trace test for the fiber of the projection $\hat p:\widehat Y\to \C^n$.
Intersecting the variety $\widehat Y$ with a linear space of codimension $n{-}1$ defined by general linear polynomials in $u_1,\dots,u_n$, 
yields a 1-dimensional subvariety of $\widehat Y$ in $\C^{n+m}\times \C^n$.
Let $\cL$ be a bilinear space defined by the product of two general affine linear polynomials $l_1(x)$ and $l_2(u)$ in the unknowns $x_1,\dots,x_n$ and $u_1,\dots,u_n$ respectively.  

If we define $\cL_t$ by $l_1(x)l_2(u)+t$, the trace is affine  linear in $t$ for all $x,u$ coordinates. Also, for  any proper nonempty subset of $\cL_t\cap\widehat Y$, the coordinate-wise sum is non-linear in $t$ for at least one of the $x,u$ coordinates.

\begin{example}
In the following example, 
we examine the previously discussed trace for the  equations ~\eqref{eq:local_system} that define the total critical variety of $X$ defined by $f=x_1x_4-x_2x_3$ and objective function $\OF_u=(x_1-u_1)^2+\cdots+(x_4-u_4)^2$.

\begin{verbatim}
i2 : R=CC[lam1,x1,x2,x3,x4,u1,u2,u3,u4,t]
    --Defines our model X.
i3 : modelEqs={det matrix{{x1,x2},{x3,x4}}}  
    --Defines the total critical variety.
i4 : critEqs=modelEqs|{   (x1-u1)-lam1*diff(x1,f),   (x2-u2)-lam1*diff(x2,f),
   (x3-u3)-lam1*diff(x3,f),   (x4-u4)-lam1*diff(x4,f)}
    --Defines a codimension 3 linear space
i5 : linearSpaceU={u1-1,u2-3,u3+u4-5} 
    --the intersection of  linearSpaceU and the total critical variety is a curve.
i6 : aCurve=linearSpaceU|critEqs
    --A linear polynomial in the x-coordinates
i7 : linear1=2*x1+3*x2+5*x3+7*x4-1  
    --A linear polynomial in the u-coordinates
i8 : linear2=u3+2*u4-7
    --a bilinear polynomial in x,u coordinates.  
i9 : Lt=linear1*linear2+t
\end{verbatim}   
Note that the ED-degree of $X$ is equal to the number of points in the
intersection of \verb+aCurve+ with the variety of \verb-linear2-. When
$t=0$, the polynomial \verb-Lt- factors. The ED-degree is
the number of points of the intersection of \verb+Lt+ with \verb+aCurve+
that are also in the variety of  \verb-linear2-.
\begin{verbatim}
i10 : G={Lt}|aCurve  --G is a zero dimensional system when t is specified.
i11 : sols1=bertiniZeroDimSolve(({t}|G),MPTYPE=>2,USEREGENERATION=>1);    --t=0
i12 : sols2=bertiniZeroDimSolve(({t+.5}|G),MPTYPE=>2,USEREGENERATION=>1); --t=-.5
i13 : sols3=bertiniZeroDimSolve(({t+1}|G),MPTYPE=>2,USEREGENERATION=>1);  --t=-1
--There are seven solutions. Although the ED-degree is 2.
i14 : #sols1
o14 = 7  
--The trace of a set of solutions is the coordinate wise sum of the solution set.
i15 : trace1=sum(sols1/ coordinates)
i16 : trace2=sum (sols2/ coordinates)
i17 : trace3=sum (sols3/ coordinates)
o17 = {2.29798, .0471255, 6.97206, 20.0172, -10.9275, 7, 21,
       --------------------------------------------------------------
       28.7778, 6.22222, -7}
\end{verbatim}
The trace moves linearly in the $x,u$ coordinates. Therefore, 
the trace test is verified, if after 
dropping the \verb+lam1+ coordinate, the $x,u$ coordinates of the following 
difference is zero. 
\begin{verbatim}
i18 : drop((trace1-trace2)-(trace2-trace3),1)

oo18 = {-3.33067e-16+1.9984e-15*ii, -1.77636e-15+6.66134e-16*ii,
       --------------------------------------------------------------
       1.77636e-14-3.27516e-15*ii, -3.55271e-15+8.32667e-16*ii,
       --------------------------------------------------------------
       6.93889e-17*ii, 3.10862e-15*ii, 1.06581e-14-3.77476e-15*ii,
       --------------------------------------------------------------
       -3.55271e-15+3.55271e-15*ii, -4.996e-16*ii}
\end{verbatim}
\end{example}


\subsection{Finding the first critical point}\label{ss:1st_critpt}

To conclude with our method, we discuss now the way we find the first critical point, so that we can run the algorithms discussed in the previous sections. The first method uses gradient descent homotopies that were introduced in~\cite{Hau13} and the idea is the following. Let $u\in \C^n$ be fixed parameter values and suppose that $x_0\in X$ is chosen at random, then most likely the vector $\nabla \OF_u(x_0)$ will not be in the linear span of the columns of the Jacobian matrix $\JacX$. Thus, not all of the $(m{+}1)$-minors of the extended Jacobian matrix~\eqref{eq:crit_conditions} are zero. We use gradient descents and homotopy continuation to track $x_0$ to a point $x^*\in X$ where these determinantal conditions are satisfied.

We assume that we know at least one point $x_0\in X$. This is not hard to achieve. 
For instance, we could have a witness set, which is the intersection of $X$ with a general $m$-dimensional affine linear space $\mathcal L$. Otherwise, we could find one point in $x_0\in X$ using cheaters homotopy as follows. Let $\hat x\in \C^n$ be chosen at random and let $\alpha := F(\hat x)\in \C^m$. Most likely, $\alpha\neq 0$, thus we use the following parameter homotopy
\begin{equation}\label{eq:cheatersHomot}
h(x,t) := F(x) - (1-t) \alpha \, ;
\end{equation}
letting $t$ run from 0 to 1, we track $\hat x$ to a point $x_0$ in $X$.

In general, the starting point $x_0$ will not satisfy the polynomial system $G(x,\lambda)$ from~\eqref{eq:local_system}. Although $F(x_0)=0$, the linear combination 
$$
\lambda_0 \nabla \OF(x;u) + \lambda_1 \nabla f_1(x) +\cdots + \lambda_m \nabla f_m(x) = K,
$$
for some nonzero vector $K\in \C^n$. We consider the following system
\begin{equation}\label{eq:local_system2}
\widehat G_K(x,\lambda) := \left[\begin{array}{c} 
F(x) \\
\lambda_0 \nabla \OF(x;u) + \lambda_1 \nabla f_1(x) +\cdots + \lambda_m \nabla f_m(x) - K
\end{array}\right].
\end{equation}
Note that for any $K\in \R^n$ there exists $\lambda \in \P^m$ such that $\widehat G(x_0, \lambda) = 0$. We need an appropriate $K$ and a homotopy that let $K\to 0$.
We start by defining $ K := \nabla \OF_u(x_0)$ and the homotopy 
\begin{equation}\label{eq:homotopy1}
H_K(x,\lambda,t) := \left[\begin{array}{c} 
F(x) \\
\lambda_0 \nabla \OF_u(x) + \lambda_1 \nabla f_1(x) +\cdots + \lambda_m \nabla f_m(x) - (1-t)\cdot K
\end{array}\right].
\end{equation}
If $\lambda = [1: 0 :\dotsc : 0]$,  the point $(x_0,\lambda)$ is a solution to~\eqref{eq:homotopy1} when $t=0$. We use homotopy continuation to find a  solution $(x^*, \lambda^*)$  to~\eqref{eq:homotopy1} when $t=1$. To guarantee that $\lambda \in \P^m$,  we work inside the following affine patch of $\P^m$:
\begin{equation}\label{eq:homotopy2}
H^a_K(x,\lambda,t) := \left[\begin{array}{c} 
F(x) \\
\lambda_0 \nabla \OF_u(x) + \lambda_1 \nabla f_1(x) +\cdots + \lambda_m \nabla f_m(x) - (1-t)\cdot K\\
\lambda_0 + a_1\lambda_1+\cdots + a_m \lambda_m - a_0
\end{array}\right],
\end{equation}
where $a_i\in \R\setminus\{0\}$ are chosen at random. For the homotopy $H^a_K$, we start at $t=0$ from the point $(x_0,a_0,0,\ldots,0)$ and track it to a solution $(x^*,\lambda_0^*,\ldots, \lambda_m^*)\in \C^{n+m+1}$  for $t=1$, where not all of the $\lambda_i^*$ are zero.
The computed $x^*$ will be a critical point and it becomes the starting point of our monodromy method.

%
\section{Implementation and Illustrating Example}
In this section, we report some of the computations we have done, including those that could be achieved with our method for the first time. 
\subsection{Planar ellipse}
We start by going step-by-step in our method to compute the ED degree of the ellipse from Example~\ref{firstExample} in the Introduction. The following code defines the equations of the ellipse. For  numerical stability, we start by normalizing the coefficients of the equation.
\begin{verbatim}
i1 : R=QQ[x1,x2,u1,u2,L1];
i2 : fModel=1/3000*(1744*x1^2-2016*x1*x2-2800*x1+1156*x2^2+2100*x2+1125)
     218  2    84        289  2   14      7     3
o2 = ---x1  - ---x1*x2 + ---x2  - --x1 + --x2 + -
     375      125        750      15     10     8
o2 : R
\end{verbatim}

\subsection{Determining the first point}
We find one critical point in the ellipse as we explained in Section 3.3, by using the cheaters homotopy to find one point in the ellipse, and the gradient descent to track it to a critical point. 
First, we compute a random point $x_r\in \C^2$ and define $b = f(x_r)$. We use the cheaters homotopy~\eqref{eq:cheatersHomot} and use the {\tt track} function from the {\tt NumericalAlgebraicGeometry package}~\cite{NAG4M2} inside {\tt Macaulay2}, which requires of a squared system. We square the system by considering a random line $L$ that passes through $x_r$. We built up this function in our library, namely \verb+pSlice+, which takes a point $p$ and a number $a$ and creates a random linear space of dimension $a$ passing through $p$.
\begin{verbatim}
i10 : xrand = flatten entries random(R^2,R^1) --random point in QQ^2
       1
o10 = {-, 1}
       2
i11 : b =sub(fModel,{x1=>xrand_0, x2=>xrand_1}) --substitue in fModel
i15 : L = pSlice(xrand,1)
o15 = {.0857144x1 + .506099x2 - .548957}
o15 : List
i16 : p2 = flatten entries random(R^2,R^1);
i17 : L2 = pSlice(p2,1);
i18 : fstart = {fModel-b} | L;
i19 : fend = {fModel} | L2;
i20 : xInit = track(fstart,fend, {sols}) -- xInit is a point in the ellipse 
i21 : xInit = coordinates first xInit
o21 = {.365254+.807261*ii, .165208+.859724*ii}
o21 : List
\end{verbatim}
We use the gradient descent homotopy described in~\eqref{eq:homotopy2} to track the random point $x_r$ to a critical point. This function is implemented in our library under the name \verb|gradDescHomot|, which takes a system, an initial point $x_0\in X$ and a parameter point $u$.
\begin{verbatim}
i23 : uPoint={0.75, -0.29}; --start data point
i24 : firstCritPt = gradDescHomot({fModel}, xInit, uPoint);
o24 = {.598568, -.0941507, .573787, -.498996}
o24 : List
i25 : xPoint = take(firstCritPt, 2);
o25 = {.598568, -.0941507}
o25 : List
i26 : lambdas = take(firstCritPt,-2)
o26 = {.573787, -.498996}
o26 : List
i27 : lagPoint = {lambdas_1/lambdas_0}
o27 = {-.869654}
o27 : List
\end{verbatim}
Thus, we found the point $x^*= (0.598568, -0.0941507)$ as a critical point, and our function also gives the values of the Lagrange multipliers, so $ \lambda^* =(0.573787, -0.498996)$ and these form a solution to the system $H(x,\lambda,1)$ defined in~\eqref{eq:homotopy1}. Therefore, in the affine patch when $\lambda_0 = 1$, the value of $\lambda_1 = -0.869654$ and the system  $H(x,\lambda,1)$ becomes the system~\eqref{eq:criticalsystem} from the introduction.

\subsection{Performing monodromy}
Now that we have one critical point $(x^*,\lambda^*)$, we compute random monodromy paths to find all the complex critical points.
\begin{verbatim}
i28 : ourStartPoint=xPoint|lagPoint
o28 - {.598567767156653, -.0941507204851288, -.869654152906697}
o28 : List
\end{verbatim}
We call Bertini, so we need to specify the directory where we want Bertini and Macaulay2 store the temporary files. We chose the variable \verb+theDir+ to store the string with the directory.
\begin{verbatim}
i41 : M=matrix{{x1-u1,x2-u2}}||matrix{{diff(x1,fModel),diff(x2,fModel)}}
o41 = | x1-u1                    x2-u2                    |
      | 436/375x1-84/125x2-14/15 -84/125x1+289/375x2+7/10 |
              2       2
o41 : Matrix R  <--- R
i42 : critEqs=flatten entries ( matrix{{1,L1}}*M)
       436         84                  14       84        289                   7
o42 = {---x1*L1 - ---x2*L1 + x1 - u1 - --L1, - ---x1*L1 + ---x2*L1 + x2 - u2 + --L1}
       375        125                  15      125        375                  10
o42 : List
i43 : Eqs={fModel}|critEqs;

i29 : makeB'InputFile(theDir,
          B'Configs=>{
          {"MPTYPE",2},
    {"PARAMETERHOMOTOPY",2},
    {"USEREGENERATION",1}},
    AVG=>{{x1,x2,L1}},
    PG=>{u1,u2},
    B'Polynomials=>Eqs);

i31 : theSolutionsViaMonodromy=b'PHMonodromyCollect(theDir,
        MonodromyStartPoints=>{ourStartPoint},--a list of points
        MonodromyStartParameters=>uPoint,
                NumberOfLoops=>100,NumSolBound=>4);
i34 : importSolutionsFile(theDir, NameSolutionsFile=>"start")      
o34 = {{.252902, -.814004, -5.38676}, {.598568, -.0941507, -.869654}, 
       {.83295, -.00662553, -2.09671}, {.844456, -.333067,-.346872}}
\end{verbatim}
Our procedure took 2 seconds and 18 random loops to find all the four critical points 
given in Example~\ref{firstExample} from the Introduction.

\subsection{Performing the trace test}
In this section, we outline code that allows us to do a trace test.
We recall that our model is defined by the equation \texttt{fModel} and the total critical variety is defined by \texttt{Eqs}.
\begin{verbatim}
i1 : R=CC[x1,x2,L1,u1,u2,t]
i2 : fModel=1/3000*(1744*x1^2-2016*x1*x2-2800*x1+1156*x2^2+2100*x2+1125)
i3 : M=matrix{{x1-u1,x2-u2}}||matrix{{diff(x1,fModel),diff(x2,fModel)}}
i4 :  critEqs=flatten entries ( matrix{{1,L1}}*M)
i5 : Eqs={fModel}|critEqs;
\end{verbatim}
We now consider a linear space defined by \texttt{sliceU} of codimension $1$.
Intersecting the $2$ dimensional total critical variety with \texttt{sliceU} produces a curve for which we will perform the trace test on. 
We intersect this curve with \texttt{Lt} which is a bilinear polynomial in the $x,u$ coordinates. This intersection consists of $6$ points.  When  $t=0$, $4$ of the points vanish on \texttt{l2} and correspond to the $4$ points we found in the previous subsection.  
Now we can use these four points as start points for a \textit{second} monodromy computation. This monodromy computation treats $t$ as the parameter and the $x's,u's,\lambda's$ as unknowns. 
\begin{verbatim}
i6 : sliceU=.3*(u1-.75)-.1*(u2+.29)
i7 : l1 = (.2*x1+.3*x2+.5)
i8 : l2 = (u1-.75)+.7*(u2+.29)
i11 : Lt=l1*l2+t
i13 : printingPrecision=200
i14 : startPointXLU={.598567767156653, -.0941507204851288, -.869654152906697,
      .75,-.29}
i15 : EqsAll={sliceU,Lt}|Eqs
i16 : makeB'InputFile(theDir,
                B'Configs=>{
                {"MPTYPE",2},
          {"PARAMETERHOMOTOPY",2},
          {"USEREGENERATION",1}},
          AVG=>{{x1,x2,L1,u1,u2}},
          PG=>{t},
          B'Polynomials=>EqsAll);
i17 :  solutionsForTraceTest=b'PHMonodromyCollect(theDir,
              MonodromyStartPoints=>{startPointXLU},--a list of points
              MonodromyStartParameters=>{0},
                      NumberOfLoops=>50,NumSolBound=>6)
\end{verbatim}
We perform a trace test to verify we have found all solutions. 
By deforming $t$ linearly from 0 to $.1\gamma$ and to $.2\gamma$, we now have three solutions sets.  Storing the traces of these solution sets as \texttt{threeTraces} we do the trace test at \texttt{i25}. The output is numerically zero for the $x$ and $u$ coordinates showing that the trace is indeed linear. Therefore we conclude that we have found all of the points.  
And because at $t=0$ there are $4$ of six solutions on $\texttt{l2}$ we conclude the ED-degree is $4$.
\begin{verbatim}                      
i18 : gamma=.0177494619790914+.60014762266504*ii
i19 : threeSolutionSets={ solutionsForTraceTest}
i20 : b'PHSequence(theDir,{{.1*gamma}})
i21 : threeSolutionSets=append(threeSolutionSets, importSolutionsFile(theDir));
i22 : b'PHSequence(theDir,{{.2*gamma}})
i23 : threeSolutionSets=append(threeSolutionSets, importSolutionsFile(theDir));
i24 : threeTraces=for i in threeSolutionSets list sum i
i25 : (threeTraces_0-threeTraces_1)-(threeTraces_1-threeTraces_2)
oo25 = {-4.44089209850063e-15+3.88578058618805e-16*ii,
        -------------------------------------------------------------
        -4.44089209850063e-15+2.54787510534094e-16*ii,
        -------------------------------------------------------------
        -1.24344978758018e-14-3.33066907387547e-15*ii,
        -------------------------------------------------------------
        4.44089209850063e-16+4.57966997657877e-16*ii,
        -------------------------------------------------------------
        1.77635683940025e-15-7.21644966006352e-16*ii}
\end{verbatim}

\section{Experimental Timings}
We end  with experimental timings  showing the monodromy method effectively computes critical points. 
Note that these timings were done in serial but can also be performed in~parallel. 
\subsection{Matrix Models in Statistics}
In \cite{HRS12}, a polynomial system to determine the critical points of likelihood equations is given for matrices with rank constraints. 
By solving this polynomial system using regeneration methods \cite{HSW11}, the ML-degrees for various $m\times n$ matrices $[p_{ij}]$ with rank at most $r$  were computed. 
Their results include detraining the following ML-degrees. 
$$
\begin{matrix}
        & (m,n) = & (3,3) & (3,4) & (3,5) & (4,4) & (4,5)  & (4,6) &\\
r=1 & &      1 &   1   &   1 & 1 &  1 &  1     & \\
r=2 & &      10 & 26 & { 58} & { 191} & { 843} & { 3119} & \\
r=3 & &        1         &1 &1 & { 191} &  { 843} &  { 3119} & \\
r=4 & &       & & &1 &1 &1 &
\end{matrix}
$$
In \cite{HRS12},  the following experimental timings for the rows labeled * in Table \ref{tab:PreprocessTime} were made.
We include our experimental timings for computing this number of critical points in bold using the monodromy method.
We see that the monodromy method performs significantly  faster in determining these critical points. We also include the number of loops that were needed to compute these critical points.

\begin{table}[htb]  \centering
  \begin{tabular}{ccccccc}
  $(m,n,r)$ & $(3,3,2)$ & $(3,4,2)$ & $(3,5,2)$ & $(4,4,2)$ & $(4,4,3)$ \\
  \hline
  Polyhedral using {\tt PHC}* & 4s & 120s & 2017s & 23843s & 1869s \\
  Regeneration using {\tt Bertini}* & 6s & 61s & 188s & 2348s & 7207s \\
  Monodromy using {\tt Bertini.m2} & \textbf{4s} &\textbf{10s}  & \textbf{79s} & \textbf{322s} & \textbf{496s} \\
  and \# of monodromy loops         & 8  & 11    & 11 & 13 & 18 \\
\end{tabular}
\caption{
Running times for preprocessing in serial.  
 Rows marked * are from \cite{HRS12}.}\label{tab:PreprocessTime}
\end{table}

\bibliographystyle{amsplain}
\bibliography{refs}

\providecommand{\bysame}{\leavevmode\hbox to3em{\hrulefill}\thinspace}
\providecommand{\MR}{\relax\ifhmode\unskip\space\fi MR }
\providecommand{\MRhref}[2]{%
  \href{http://www.ams.org/mathscinet-getitem?mr=#1}{#2}
}
\providecommand{\href}[2]{#2}
\begin{thebibliography}{10}

\bibitem{AAGL14}
D.~{Agostini}, D.~{Alberelli}, F.~{Grande}, and P.~{Lella}, \emph{{The maximum
  likelihood degree of {F}ermat hypersurfaces}}, preprint arXiv:1404.5745
  (2014).

\bibitem{BBHSW14}
D.~J. Bates, D.~A. Brake, J.~D. Hauenstein, A.~J. Sommese, and C.~W. Wampler,
  \emph{Homotopies for connected components of algebraic sets with application
  to computing critical sets},  (2014).

\bibitem{Bertini4M2}
D.~J. Bates, E.~Gross, A.~Leykin, and J.~I. Rodriguez, \emph{{Bertini for
  Macaulay2}}, preprint arXiv:1310.3297 (2013).

\bibitem{Bertini}
D.~J. Bates, J.~D. Hauenstein, A.~J. Sommese, and C.~W. Wampler,
  \emph{{Bertini: Software for Numerical Algebraic Geometry}}, Available at
  bertini.nd.edu with permanent doi: dx.doi.org/10.7274/R0H41PB5.

\bibitem{BHSW13}
\bysame, \emph{Numerically solving polynomial systems with {B}ertini},
  Software, Environments, and Tools, vol.~25, Society for Industrial and
  Applied Mathematics (SIAM), Philadelphia, PA, 2013. \MR{3155500}

\bibitem{CHKS06}
F.~Catanese, S.~Ho{\c{s}}ten, A.~Khetan, and B.~Sturmfels, \emph{The maximum
  likelihood degree}, American Journal of Mathematics \textbf{128} (2006),
  no.~3, 671--697.

\bibitem{DHOST13}
J.~{Draisma}, E.~{Horobet}, G.~{Ottaviani}, B.~{Sturmfels}, and R.~R. {Thomas},
  \emph{{The {E}uclidean distance degree of an algebraic variety}}, preprint
  arXiv:1309.0049 (2013).

\bibitem{DR14}
J.~Draisma and J.~I. Rodriguez, \emph{Maximum likelihood duality for
  determinantal varieties}, International Mathematics Research Notices
  \textbf{2014} (2014), no.~20, 5648--5666.

\bibitem{DLR15}
D.~{Drusvyatskiy}, H.-L. {Lee}, and R.~R. {Thomas}, \emph{{Counting real
  critical points of the distance to spectral matrix sets}}, preprint
  arXiv:1502.02074 (2015).

\bibitem{FauSafSpa12}
J.-C. Faug\`ere, M.~{Safey El Din}, and P.-J. Spaenlehauer, \emph{Critical
  points and {G}r\"obner bases: the unmixed case}, Proceedings of the 2012
  International Symposium on Symbolic and Algebraic Computation (ISSAC 2012),
  2012, pp.~162--169.

\bibitem{M2}
D.~R. Grayson and M.~E. Stillman, \emph{Macaulay2, a software system for
  research in algebraic geometry}, Available at
  http://www.math.uiuc.edu/Macaulay2/.

\bibitem{GR13}
E.~Gross and J.~I. Rodriguez, \emph{Maximum likelihood geometry in the presence
  of data zeros}, Proceedings of the 39th International Symposium on Symbolic
  and Algebraic Computation (New York, NY, USA), ISSAC '14, ACM, 2014,
  pp.~232--239.

\bibitem{Hau13}
J.~D. Hauenstein, \emph{Numerically {C}omputing {R}eal {P}oints on {A}lgebraic
  {S}ets}, Acta Applicandae Mathematicae \textbf{125} (2013), no.~1, 105--119.

\bibitem{HR15}
J.~D. Hauenstein and J.~I. Rodriguez, \emph{Numerical irreducible decomposition
  for multiprojective varieties}, In preparation.

\bibitem{HRS12}
J.~D. {Hauenstein}, J.~I. {Rodriguez}, and B.~{Sturmfels}, \emph{{Maximum
  Likelihood for Matrices with Rank Constraints}}, preprint arXiv:1210.0198
  (2012).

\bibitem{HSW11}
J.~D. Hauenstein, A.~J. Sommese, and C.~W. Wampler, \emph{Regenerative cascade
  homotopies for solving polynomial systems}, Applied Mathematics and
  Computation \textbf{218} (2011), no.~4, 1240 -- 1246.

\bibitem{HKS05}
S.~Ho{\c{s}}ten, A.~Khetan, and B.~Sturmfels, \emph{Solving the likelihood
  equations}, Found. Comput. Math. \textbf{5} (2005), no.~4, 389--407.
  \MR{2189544}

\bibitem{Huh13}
J.~Huh, \emph{The maximum likelihood degree of a very affine variety}, Compos.
  Math. \textbf{149} (2013), no.~8, 1245--1266. \MR{3103064}

\bibitem{HS13}
J.~{Huh} and B.~{Sturmfels}, \emph{{Likelihood Geometry}}, preprint
  arXiv:1305.7462 (2013).

\bibitem{John48}
F.~John, \emph{Extremum problems with inequalities as subsidiary conditions},
  Studies and {E}ssays {P}resented to {R}. {C}ourant on his 60th {B}irthday,
  {J}anuary 8, 1948, Interscience Publishers, Inc., New York, N. Y., 1948,
  pp.~187--204. \MR{0030135 (10,719b)}

\bibitem{Lee14}
H.~{Lee}, \emph{{The Euclidean Distance Degree of {F}ermat Hypersurfaces}},
  preprint arXiv:1409.0684 (2014).

\bibitem{NAG4M2}
A.~Leykin, \emph{Numerical algebraic geometry}, Journal of Software for Algebra
  and Geometry \textbf{3} (2011), no.~1, 5--10.

\bibitem{NR09}
J.~Nie and K.~Ranestad, \emph{Algebraic degree of polynomial optimization},
  SIAM J. Optim. \textbf{20} (2009), no.~1, 485--502. \MR{2507133
  (2010d:90105)}

\bibitem{OttSpaStu14}
G.~Ottaviani, P.-J. Spaenlehauer, and B.~Sturmfels, \emph{Exact solutions in
  structured low-rank approximation}, SIAM Journal on Matrix Analysis and
  Applications \textbf{35} (2014), no.~4, 1521--1542.

\bibitem{Pie15}
R.~Piene, \emph{Polar varieties revisited}, Computer Algebra and Polynomials
  (J.~Gutierrez, J.~Schicho, and M.~Weimann, eds.), Lecture Notes in Computer
  Science, Springer International Publishing, 2015, pp.~139--150.

\bibitem{Ran12}
K.~Ranestad, \emph{Algebraic degree in semidefinite and polynomial
  optimization}, Handbook on semidefinite, conic and polynomial optimization,
  Internat. Ser. Oper. Res. Management Sci., vol. 166, Springer, New York,
  2012, pp.~61--75. \MR{2894691}

\bibitem{Rod14}
J.~I. Rodriguez, \emph{Maximum likelihood for dual varieties}, Proceedings of
  the 2014 Symposium on Symbolic-Numeric Computation (New York, NY, USA), SNC
  '14, ACM, 2014, pp.~43--49.

\bibitem{RT15}
J.~I. {Rodriguez} and X.~{Tang}, \emph{{Data-Discriminants of Likelihood
  Equations}}, preprint arXiv:1501.00334 (2015).

\bibitem{SVW2002}
A.~J. Sommese, J.~Verschelde, and C.~W. Wampler, \emph{Symmetric functions
  applied to decomposing solution sets of polynomial systems}, SIAM J. Numer.
  Anal. \textbf{40} (2002), no.~6, 2026--2046 (2003). \MR{1974173
  (2004m:65069)}

\bibitem{SW05}
A.~J. Sommese and C.~W. {Wampler, II}, \emph{The numerical solution of systems
  of polynomials}, World Scientific Publishing Co. Pte. Ltd., Hackensack, NJ,
  2005, Arising in engineering and science. \MR{2160078 (2007a:14065)}

\bibitem{Spa14}
P.-J. Spaenlehauer, \emph{On the complexity of computing critical points with
  {G}r\"obner bases}, SIAM J. Optim. \textbf{24} (2014), no.~3, 1382--1401.
  \MR{3248045}

\bibitem{WR13}
W.~Wu and G.~Reid, \emph{Finding points on real solution components and
  applications to differential polynomial systems}, I{SSAC}
  2013---{P}roceedings of the 38th {I}nternational {S}ymposium on {S}ymbolic
  and {A}lgebraic {C}omputation, ACM, New York, 2013, pp.~339--346.
  \MR{3206376}

\end{thebibliography}

\end{document}